\documentclass[12pt,a4paper]{amsart}
\usepackage{amssymb,amsmath}

\headsep=1cm \numberwithin{equation}{section}
\newtheorem{theorem}{Theorem}[section]
\newtheorem{lemma}[theorem]{Lemma}
\newtheorem{proposition}[theorem]{Proposition}
\newtheorem{corollary}[theorem]{Corollary}

\theoremstyle{definition}
\newtheorem{definition}[theorem]{Definition}

\newtheorem{remark}[theorem]{Remark}

\newcommand\Supp{\operatorname{Supp}}
\newcommand\Ass{\operatorname{Ass}}

\newcommand\Ann{\operatorname{Ann}}
\newcommand\Spec{\operatorname{Spec}}
\newcommand\Rad{\operatorname{Rad}}

\begin{document}

\title[strongly irreducible submodules]{Some characterizations of strongly irreducible submodules in arithmetical and Noetherian modules}%
\author{Reza Naghipour$^*$ and Monireh Sedghi}%
\address{Department of mathematics, university of tabriz, tabriz, iran, and School of Mathematics,
 Institute for Research in Fundamental Sciences (IPM), P.O. Box: 19395-5746, Tehran, Iran.}%
\email{naghipour@ipm.ir (R. NAGHIPOUR)}%
\address{Department of Mathematics, Azarbaijan Shahid Madani University, Tabriz, Iran.}
\email{m\_sedghi@tabrizu.ac.ir (M. SEDGHI)}%
\thanks{
$^*$Corresponding author: e-mail: {\it naghipour@ipm.ir} (Reza Naghipour).}%

\begin{abstract}
The purpose of the present paper is to prove some properties of the
strongly irreducible submodules in the arithmetical and Noetherian
modules over a commutative ring. The relationship among the families
of strongly irreducible submodules,  irreducible submodules, prime
submodules and primal submodules is proved. Also, several new
characterizations of the arithmetical modules are given. In the case
when $R$ is Noetherian and $M$ is finitely generated, several
characterizations of strongly irreducible submodules are included.
Among other things, it is shown that when $N$ is a submodule of $M$
such that $N:_RM$ is not a prime ideal, then $N$ is strongly
irreducible if and only if there exist submodule $L$ of $M$ and
prime ideal $\frak p$ of $R$ such that $N$ is $\frak p$-primary,
$N\subsetneqq L\subseteq \frak pM$ and for all submodules $K$ of $M$
either $K\subseteq N$ or $L_{\frak p}\subseteq K_{\frak p}$. In
addition, we show that a submodule $N$  of $M$ is strongly
irreducible if and only if $N$ is primary, $M_{\frak p}$ is
arithmetical and $N=(\frak pM)^{(n)}$ for some integer $n>1$, where
$\frak p=\Rad(N:_RM)$ with $\frak p\not\in \Ass_RR/\Ann_R(M)$ and
$\frak pM\nsubseteq N$. As a consequence we deduce that if $R$ is
integral domain and $M$ is torsion-free, then there exists a
strongly irreducible submodule $N$ of $M$ such that $N:_RM$ is not
prime ideal if and only if there is a prime ideal $\frak p$ of $R$
with $\frak pM\nsubseteq N$ and $M_{\frak p}$ is an arithmetical
$R_{\frak p}$-module.
\end{abstract}
\maketitle
\section {Introduction}
Let $R$ be a commutative ring with non-zero identity and let $M$ be
an arbitrary $R$-module. We say that a submodule $N$ of $M$ is a
{\it distributive submodule} if for all submodules $K, L$ of $M$,
the following equivalent conditions are satisfied:

${\rm (i)}$ $(K+L)\cap N=(K\cap N)+(L\cap N);$

${\rm (ii)}$ $(K\cap L)+ N=(K+ N)\cap(L+ N).$

 Also, $M$ is said to be {\it distributive module} if every submodule of $M$ is a distributive
submodule.

We say that $N$ is an {\it irreducible submodule} if $N$ is not the
intersection of two submodules of $M$ that properly contain it. It
is easy to see that if an irreducible submodule $N$ is distributive,
then for all submodules $K, L$ of $M$ the condition $K\cap
L\subseteq N$ implies that either $K\subseteq N$ or $L\subseteq N$.
These considerations motivated us to define a submodule $N$ of an
$R$-module $M$ to be {\it strongly irreducible}, if for all
submodules $K, L$ of $M$, the condition $K\cap L\subseteq N$ implies
that either $K\subseteq N$ or $L\subseteq N$.

The purpose of the present article is to introduce and examine some
properties of strongly irreducible submodules. In particular, we
relate the notions of strongly irreducible submodules, irreducible
submodules, primal submodules and prime submodules of an $R$-module
$M$. Also a characterization of strongly irreducible submodules in
an arithmetical module is given. Specifically, we show that a
submodule $N$ of an arithmetical $R$-module $M$ is strongly
irreducible if and only if $N$ is a primal submodule if and only if
$N$ is an irreducible submodule. If moreover $R$ is assumed to be
Noetherian ring and $M$ finitely generated, then $N$ is strongly
irreducible if and only if $N$ is primary, $M_{\frak p}$ is
arithmetical and $N=(\frak pM)^{(n)}$ for some integer $n>1$, where
$\frak p={\rm Rad}(N:_RM)$ with $\frak p\not\in \Ass_R(R/\Ann_R(M))$
and $\frak pM\nsubseteq N$.

We recall that an $R$-module $M$ is said to be {\it arithmetical
module} if $M_{\frak m}$ is an uniserial module over $R_{\frak m}$,
for each maximal ideal $\frak m$ of $R$, i.e., the submodules of
$M_{\frak m}$ are linearly ordered with respect to inclusion.

A brief summary of the contents of this article will now be given.
Let $R$ be a commutative ring and let $M$ be an arbitrary
$R$-module. In Section 2, the notion of strongly irreducible
submodules are introduced, and some properties of them are
considered. A typical result in this direction, which is a
generalization of the main result of Heinzer et al. (see
\cite[Theorem 2.6]{HRR}) to strongly irreducible submodules,  is the
following:
\begin{theorem}
Let $(R, \frak m)$ be a quasi-local ring and let $M$ be an
$R$-module. Suppose that $N$ is a strongly irreducible submodule of
$M$ such that $N\neq N:_M\frak m$. Then the following conditions
hold:

${\rm (i)}$ the submodule $N:_M\frak m$ of $M$ is cyclic,

${\rm (ii)}$ $N=\frak m(N:_M\frak m)$,

${\rm (iii)}$ for each submodule $K$ of $M$ either $K\subseteq N$ or
$ N:_M\frak m\subseteq K$.
\end{theorem}

The proof of Theorem 1.1 is given in Theorem 2.11. Using this, it is
shown that if $N$ is a strongly irreducible $\frak m$-primary
submodule of a finitely generated module $M$ over a local Noetherian
ring $(R, \frak m)$ with $\frak mM\neq N$, then
\begin{center}
$N=\bigcup \{K|\, K\, \text{is a submodule of}\,\, M\,
\text{and}\,\, K\subsetneqq N:_M\frak m\},$
\end{center}
 and
\begin{center}
$N:_M\frak m=\bigcap \{L|\, L\, \text{is a submodule of}\,\, M\,
\text{and}\,\, N\subsetneqq L\}.$
\end{center}

In Section 3, we give some characterizations of arithmetical
 and distributive modules.  We also establish two characterizations whenever a
submodule of an arithmetical module over a commutative ring is
strongly irreducible. More precisely, we shall prove:
\begin{theorem}
Let $M$ be an arithmetical $R$-module and let $N$ be a submodule of
$M$. Then the following statements are equivalent:

${\rm (i)}$ $N$ is strongly irreducible.

${\rm (ii)}$ $N$ is primal.

${\rm (iii)}$ $N$ is irreducible.
\end{theorem}

Finally, the main results of the Section 4, provide several
characterizations for a finitely generated module $M$ over a
Noetherian ring $R$ to have a strongly irreducible submodule. In
this section, among other things, we shall show that:
\begin{theorem}
Let $M$ be a finitely generated module over a Noetherian ring $R$
and let $N$ be a submodule of $M$. Then, $N$ is strongly irreducible
if and only if $N$ is primary, $M_{\frak p}$ is an arithmetical
$R_{\frak p}$-module and $N=(\frak pM)^{(n)}$ for some integer
$n>1$, where $\frak p={\rm Rad}(N:_RM)$ such that $\frak
pM\nsubseteq N$ and $\frak p \not\in \Ass_R(R/\Ann_R(M))$.
\end{theorem}

The proof of Theorem 1.3 is given in Theorem 4.7. One of our tools
for proving Theorem 1.3 is the following:
\begin{proposition}
Let $M$ be a non-zero finitely generated module over a local
(Noetherian) ring $(R, \frak m)$ such that ${\rm
deth}_R(R/\Ann_R(M))>0$ and $\Gamma_{\frak m}(M)\neq0$. Let $N$ be a
strongly irreducible submodule of $M$ such that ${\rm
Rad}(N:_RM)=\frak m$. Then $\Gamma_{\frak m}(M)\nsubseteq N$.
\end{proposition}
Pursuing this point of view further we derive the following
consequence of Theorem 1.3.

\begin{corollary}
Let $M$ be a non-zero finitely generated torsion-free module over a
Noetherian integral domain $R$. Then there exists a strongly
irreducible submodule $N$ of $M$ such that $N:_RM$ is not prime
ideal of $R$ if and only if there is a prime ideal $\frak p$ of $R$
with $\frak pM\nsubseteq N$ and $M_{\frak p}$ is an arithmetical
$R_{\frak p}$-module.
\end{corollary}

Finally, using Theorem 1.1 we obtain the following proposition which
gives us a characterization of a strongly irreducible submodule in a
multiplication module over a commutative Noetherian ring.

\begin{proposition}
Let $M$ be a multiplication module over a Noetherian ring $R$, and
suppose that $N$ is a proper submodule of $M$ such that the ideal
$N:_RM$ of $R$ is not prime. Then $N$ is strongly irreducible if and
only if there exist submodule $L$ of $M$ and prime ideal $\frak p$
of $R$ such that $N$ is $\frak p$-primary, $N\subsetneqq L\subseteq
\frak pM$ and for all submodules $K$ of $M$ either $K\subseteq N$ or
$L_{\frak p}\subseteq K_{\frak p}$.
\end{proposition}
Throughout this paper, $R$ will always be a commutative ring with
non-zero identity and $\frak a$  will be an ideal of $R$. For each
$R$-module $M$ and for any submodule $N$ of $M$, the submodules
$\bigcup _{s\in S}(N:_Ms)$ and $\bigcup _{n\geq 0}(0:_M \frak a^n)$
of $M$ are denoted by $S(N)$ and $\Gamma_{\frak a}(M)$ respectively,
where $S$ is a multiplicatively closed subset of $R$. In the case
$S=R\backslash \bigcup\{\frak p\in {\rm mAss}_RM/\frak aM\}$, for
any integer $m\geq 1$,  the submodule $S(\frak a^mM)$ is denoted by
$(\frak aM)^{(m)}$. The {\it radical} of $\frak a$, denoted by
$\Rad(\frak a)$, is defined to be the set $\{r \in R: r^n \in \frak
a \,\,\text{for some}\,\, n\in \mathbb{N}\}$. Finally, for any
$R$-module $L$, we shall use $Z_R(L)$ (resp. ${\rm mAss}_RL$) to
denote the set of zerodivisors on $L$ in $R$ (resp. the set of
minimal elements of $\Ass_RL$).

Let $R$ be a Noetherian ring and let $G$ be a non-zero finitely
generated $R$-module. For $\frak p \in \Supp(G)$, the $G$-height of
$\frak p$, denoted by ${\rm ht}_G\frak p$, is defined to be the
supremum of lengths of chains of prime ideals of $\Supp(G)$
terminating with $\frak p$. We shall say that an ideal $\frak a$ of
$R$ is $G$-proper if $G/\frak aG\neq 0$, and, when this is the case,
we define the $G$-height of $\frak a$ (written ${\rm ht}_G \frak a$)
to be $\inf\{{\rm ht}_G\frak p:\, \frak p\in \Supp(G)\cap V(\frak
a)\}$, where $V(\frak a)$ denotes the set $\{\frak p\in \Spec(R):\,
\frak p\supseteq \frak a\}$. Also, if $(R, \frak m)$ is local then
we use ${\rm depth}_R(G)$ to denote the maximum length of all
$G$-sequences contained in $\frak m$.

 A proper submodule $P$ of an $R$-module $M$ is said to be {\it
prime} if whenever $rx\in P$ for $r\in R$ and $x\in M$, then $x\in
P$ or $r\in (P:_RM)$.  (For more information about prime submodules,
see \cite{MM}, \cite{THS}).

An $R$-module $M$ is said to be a {\it multiplication module} if
every submodule of $M$ is of the form $\frak bM$ for some ideal
$\frak b$ of $R$. One can easily check that $M$ is a multiplication
module if and only if, for all submodules $K$ of $M$, we have
$K=(K:_RM)M$.

For any unexplained notation and terminology we refer the reader to
\cite{BH} or \cite{Mat}.

\section{Strongly irreducible submodules}
Throughout this section, $R$ will denote a commutative ring (with
identity). The purpose of this section is to introduce the concept
of strongly irreducible submodules. Several properties of them are
considered. The main goals of this section are Theorem 2.11 and
Corollary 2.12. We begin with

\begin{definition}
 Let $M$ be an $R$-module and let $N$ be a submodule of $M$. We say that $N$ is a {\it
strongly irreducible submodule} of $M$ if  for every two submodules
$L$ and $K$ of $M$, the inclusion $L\cap K\subseteq N$ implies that
either $L\subseteq N$ or $K\subseteq N$.
\end{definition}

The first lemma shows that every strongly irreducible submodule  in
a Noetherian $R$-module $M$ is primary.
\begin{lemma}
Let $M$ be an $R$-module and let $N$ be a strongly irreducible
submodule of $M$. Then $N$ is irreducible. In particular, if $M$ is
Noetherian, then $N$ is a primary submodule of $M$.
\end{lemma}

{\it Proof}. Suppose that $N$ is a strongly irreducible submodule of
$M$. Let $L$ and $K$ be submodules of $M$ such that $N=L\cap K$.
Then $L\cap K\subseteq N$. Since $N$ is strongly submodule,  it
follows that either $L\subseteq N$ or $K\subseteq N$, and so either
$L= N$ or $K=N$. Hence $N$ is irreducible. Now, the second part
follows from
the proof of \cite[Theorem 6.8]{Mat}.   \qed\\

In the next lemma we observe that every prime submodule in a
multiplication $R$-module $M$ is strongly irreducible.
\begin{lemma}
Let $M$ be a multiplication $R$-module and let $N$ be a prime
submodule of $M$. Then $N$ is strongly irreducible submodule.
\end{lemma}
{\it Proof}.  Let $L$ and $K$ be submodules of $M$ such that $L\cap
K\subseteq N$. Then $(L\cap K):_RM\subseteq (N:_RM)$, and so
$(L:_RM)\cap (K:_RM)\subseteq (N:_RM)$. As $N:_RM$ is a prime ideal
of $R$, it follows that either $L:_RM\subseteq N:_RM$ or
$K:_RM\subseteq N:_RM$. Hence, either $(L:_RM)M\subseteq (N:_RM)M$
or $(K:_RM)M\subseteq (N:_RM)M$. Now, since $M$ is a multiplication
module, it follows that either $L\subseteq N$ or $K\subseteq N$, as
required.  \qed\\

\begin{remark}
Before bringing the next result we fix a notation, which is employed
by P. Schenzel in \cite{Sc} in the case $M=R$. Let $S$ be a
multiplicatively closed subset of $R$. For a submodule $K$ of $M$,
we use $S(K)$ to denote the submodule $\bigcup_{s\in S}(K:_Ms)$.

In particular, for any ideal $\frak a$ of $R$, if $S=R\setminus
\bigcup \{\frak p\in {\rm mAss}_RM/\frak aM \}$ then for every $n\in
\mathbb{N}$; $S(\frak a^nM)$ is denoted by $(\frak aM)^{(n)}$.
\end{remark}

\begin{lemma}
Let $M$ be an $R$-module and let $N$ be a submodule of $M$. Let $S$
be a multiplicatively closed subset of $R$ such that the submodule
$S^{-1}N$ in $S^{-1}M$ is strongly irreducible. Then $S(N)$ is also
a strongly irreducible submodule of $M$.
\end{lemma}
{\it Proof}. Let $L$ and $K$ be two submodules of $M$ such that
$L\cap K\subseteq S(N)$. Then $S^{-1}(L\cap K)\subseteq
S^{-1}(S(N))$, and so it is easy to see that $S^{-1}(L)\cap
S^{-1}(K)\subseteq S^{-1}(N)$. Now, the hypothesis on $S^{-1}N$
implies that either $S^{-1}(L)\subseteq S^{-1}(N)$ or
$S^{-1}(K)\subseteq S^{-1}(N)$. Hence, either $S^{-1}(L)\cap
M\subseteq S^{-1}(N)\cap M$ or $S^{-1}(K)\cap M\subseteq
S^{-1}(N)\cap M$; and so either $L\subseteq S(N)$ or $K\subseteq
S(N)$, as required.  \qed\\

The following proposition shows that the notion of strongly
irreducible submodule behaves well under localization.
\begin{proposition}
Let $M$ be an $R$-module and let $N$ be a strongly irreducible
primary submodule of $M$. Let $S$ be a multiplicatively closed
subset of $R$ such that $\Rad(N:_RM)\cap S=\emptyset$. Then
$S^{-1}N$ is a strongly irreducible submodule of $S^{-1}M$.
\end{proposition}

{\it Proof}. Suppose that $\Omega_1$ and $\Omega_2$ are two
submodules of $S^{-1}M$ such that $\Omega_1\cap \Omega_2\subseteq
S^{-1}N$. Then, in view of \cite[Ex. 9.11]{Sh}, there exist
submodules $K$ and $L$ of $M$ such that $\Omega_1=S^{-1}K$ and
$\Omega_2=S^{-1}L$. Hence $S^{-1}K\cap S^{-1}L\subseteq S^{-1}N$,
and so $S(K)\cap S(L)\subseteq S(N)$. Since $N$ is primary and
$\Rad(N:_RM)\cap S=\emptyset$, one easily sees that $S(N)=N$. Whence
$S(K)\cap S(L)\subseteq N$. Now, as $N$ is strongly irreducible, it
follows that either $S(K)\subseteq N$ or $S(L)\subseteq N$.
Therefore,  either $S^{-1}(S(K))\subseteq S^{-1}N$ or
$S^{-1}(S(L))\subseteq S^{-1}N$, and so either $S^{-1}K\subseteq
S^{-1}N$ or $S^{-1}L\subseteq S^{-1}N$, as required.  \qed\\

\begin{lemma}
Let $\frak p$ be a prime ideal of $R$ and let $M$ be an $R$-module.
Suppose that $N$ is a $\frak p$-primary submodule of $M$ such that
$N_{\frak p}$ is a strongly irreducible submodule of $M_{\frak p}$.
Then $N$ is a strongly irreducible submodule of $M$.
\end{lemma}

{\it Proof}. Let $L$ and $K$ be two submodules of $M$ such that
$L\cap K\subseteq N$. Then $L_{\frak p}\cap K_{\frak p}\subseteq
N_{\frak p}$. Since $N_{\frak p}$ is strongly irreducible submodule,
it follows that either $L_{\frak p}\subseteq N_{\frak p}$ or
$K_{\frak p}\subseteq N_{\frak p}$. Now, as $N$ is $\frak p$-primary
submodule, it readily follows that either  $L\subseteq N$ or
$K\subseteq N$, as required.   \qed\\

The next lemma investigates how the strongly irreducible property
behaves under the faithfully flat extensions.
\begin{lemma}
Assume that $T$ is a commutative ring which is a faithfully flat
$R$-algebra. Let $M$ be an $R$-module and assume that $N$ is a
submodule of $M$ such that $N\otimes_RT$ is a strongly irreducible
submodule of $M\otimes_RT$. Then $N$ is a strongly irreducible
submodule of $M$.
\end{lemma}

{\it Proof}. Suppose that $L$ and $K$ are two submodules of $M$ such
that $L\cap K\subseteq N$. Then, in view of \cite[Theorem 7.4]{Mat},
$(L\otimes_RT)\cap (K\otimes_RT)\subseteq N\otimes_RT$. Now, as
$N\otimes_RT$ is a strongly irreducible submodule of $M\otimes_RT$,
it follows that either $L\otimes_RT \subseteq N\otimes_RT$ or
$K\otimes_RT \subseteq N\otimes_RT$.  Now, if $L\otimes_RT \subseteq
N\otimes_RT$, then we have
$$(L+N)/N\otimes_RT=((L\otimes_RT)+(N\otimes_RT))/(N\otimes_RT)=0,$$
and so by the faithfully flatness of $T$ over $R$ we have
$L\subseteq N$. A similar argument also shows that if $K\otimes_RT
\subseteq N\otimes_RT$, then $K\subseteq N$. This completes the
proof.  \qed\\

\begin{lemma}
Let $M$ be an $R$-module and let $U$ be a submodule of $M$. Assume
that $N$ is a strongly irreducible submodule of $M$ containing $U$.
Then $N/U$ is a strongly irreducible submodules of $M/U$.
\end{lemma}

{\it Proof}. Suppose that $L$ and $K$ are two submodules of $M$ such
that $U\subseteq L\cap K$ and $K/U\cap  L/U \subseteq N/U$. Then
$K\cap L\subseteq N$, and so as $N$ is strongly irreducible it
follows that either $K\subseteq N$ or $L\subseteq N$. Hence either
$K/U\subseteq N/U$ or $L/U\subseteq N/U$, as required.  \qed\\

\begin{proposition}
Let $M$ be an $R$-module and let $N$ be a submodule of $M$. Then $N$
is strongly irreducible submodule if and only if for all cyclic
submodules $L$ and $K$ of $M$ the condition $K\cap L\subseteq N$
implies that either $K\subseteq N$ or $L\subseteq N$.
\end{proposition}
{\it Proof}. One direction is clear. To prove the converse, suppose
that $T$ and $S$ are two submodules of $M$ such that $T\cap
S\subseteq N$ and $T\nsubseteq N$. Then there exists $y\in T$ such
that $y\not\in N$. Now, for all $x\in S$, we have $Rx\cap
Ry\subseteq N$. According to hypothesis $Rx \subseteq N$, and so
$S\subseteq N$, as required.  \qed\\

We are now ready to state and prove the main theorem of this section
which provides some properties of a strongly irreducible submodule
$N$ of an arbitrary module $M$ over a quasi-local ring $(R, \frak
m)$ which is properly contained in $N:_M\frak m.$

 It is easy to see that if $R$ is Noetherian and $\Rad(N:_RM)=\frak m$, then $N$ is
properly contained in $N:_M\frak m.$ This result plays an important
role in the Section 4 where we restrict attention to the case $R$ is
a Noetherian ring.

\begin{theorem}
Let $(R, \frak m)$ be a quasi-local ring and let $M$ be an
$R$-module. Let $N$ be a strongly irreducible submodule of $M$ such
that $N\neq N:_M\frak m$. Then

${\rm (i)}$ the submodule $N:_M\frak m$ of $M$ is cyclic,

${\rm (ii)}$ $N=\frak m(N:_M\frak m)$,

${\rm (iii)}$ for each submodule $K$ of $M$ either $K\subseteq N$ or
$ N:_M\frak m\subseteq K$.
\end{theorem}
{\it Proof}. In order to show $\rm(i)$, in view of the hypothesis
$N\neq N:_M\frak m$, there exists an element $x\in N:_M\frak m$ such
that $x\not\in N$. It is enough for us to show that $N:_M\frak
m=Rx$. To do this, let $y$ be an arbitrary element of $(N:_M\frak
m)\backslash Rx$. We claim that $Rx\cap Ry\subseteq N$. To this end,
set $z\in Rx\cap Ry$. Then there exist elements $a, b\in R$ such
that $z=ax=by$. Now, if $b$ is a unit in $R$, then $y\in Rx$, which
is a contrdiction. Thus we may assume that $b$ is not unit. Then
$b\in \frak m$, and so it follows from $y\in N:_M\frak m$ that
$by\in N$. That is $z\in N$, and hence $Rx\cap Ry\subseteq N$. Now,
since $N$ is strongly irreducible it follows from $x\not\in N$ that
$Ry\subseteq N$, i.e., $y\in N$. Therefore it follows that
$N:_M\frak m=N\cup Rx$, and so $N:_M\frak m=N$ or $N:_M\frak m=Rx$.
Consequently, as $N:_M\frak m \neq N$, it yields that $N:_M\frak
m=Rx$, and so the submodule $ N:_M\frak m$ of $M$ is cyclic.

To prove (ii), in view of (i), we have $N:_M\frak m=Rx$, where $x\in
(N:_M\frak m)\backslash N$. Hence $\frak mx\subseteq N$, and so
$\frak m \subseteq N:_Rx$. As $x\not\in N$, it yields that $\frak m
=N:_Rx$. Thus $\frak m(N:_M\frak m)=(N:_Rx)x$. So it is enough for
us to show that $N\subseteq (N:_Rx)x$. To do this, let $y\in N$.
Then, as $N\subsetneqq (N:_M\frak m)=Rx$, it follows that $y=rx$,
for some $r\in R$. Since $x\not\in N$, it yields that $r$ is not
unit and so $r\in \frak m$. Hence $y\in \frak mx$, and thus $y\in
(N:_Rx)x$, as required.

Finally, in order to prove (iii) suppose that $K$ is an arbitrary
submodule of $M$ such that $(N:_M\frak m)\nsubseteq K$. We have to
show that $K\subseteq N$. To do this, let $y\in K$. Since by part
(i) $N:_M\frak m=Rx$, for some $x\in (N:_M\frak m)\backslash N$, and
 $(N:_M\frak m)\nsubseteq K$, it follows that $x\not\in K$.
Moreover, $Rx\cap Ry\subseteq N$. Because, if $w\in Rx\cap Ry$, then
there are elements $r, s\in R$ such that $w=rx=sy$. Since $y\in K$
and $x\not\in K$, it follows that $r$ is not unit, and so $r\in
\frak m$. Hence $w\in \frak mx$, i.e.,  $w\in N$.  Now, as $N$ is
strongly irreducible we deduce that either $Rx\subseteq N$ or
$Ry\subseteq N$. As $x\not\in N$, it yields that $y\in N$, and hence
$K\subseteq N$, as required.  \qed\\

Before bringing the final result of this section we recall that a
proper submodule $N$ of an $R$-module $M$ is said to be {\it
sheltered} if the set of submodules of $M$ strictly containing $N$
has a smallest member $S$, called the {\it shelter} of $N$ (see
\cite[Exercise 18, p. 238]{Bo}).

\begin{corollary}
Let $(R, \frak m)$ be a local (Noetherian) ring and let $M$ be a
finitely generated $R$-module. Suppose $N$ is a strongly irreducible
submodule of $M$ such that $\Rad(N:_RM)=\frak m$. Then $N$ is
sheltered and its shelter is $N:_M\frak m.$
\end{corollary}
{\it Proof}. Since $N$ is a strongly irreducible and
$\Rad(N:_RM)=\frak m$, it follows from Lemma 2.2 that $N$ is an
$\frak m$-primary submodule of $M$. Hence it is easy to see that
$N\neq N:_M\frak m$. Now, the assertion follows from Theorem 2.11.
\qed

\section{Strongly irreducible submodule in Arithmetical modules}

The first main result of this section, gives us several
characterizations of an arithmetical module over a commutative ring.
Before stating that theorem, let us recall that an $R$-module $M$ is
called an {\it arithmetical module} if $M_{\frak m}$ is uniserial
$R_{\frak m}$-module for each maximal ideal $\frak m$ of $R$, i.e.,
the submodules of $M_{\frak m}$ are linearly ordered with respect to
inclusion. It is clear that every primary submodule of an
arithmetical module is strongly irreducible.

\begin{theorem}
Let $M$ be an $R$-module. Then the following statements are
equivalent:

${\rm (i)}$ $M$ is an arithmetical module.

${\rm (ii)}$ $M$ is a distributive module.

${\rm (iii)}$ $(K+L):_RN=(K:_RN)+(L:_RN)$ for all submodules $K, L,
N$ of $M$ with $N$ is finitely generated.

${\rm (iv)}$ $K:_R(L\cap N)=(K:_RN)+(L:_RN)$ for all submodules $K,
L, N$ of $M$ with $L$ and $N$ are finitely generated.

${\rm (v)}$ Every finitely generated submodule of $M$ is a
multiplication module.
\end{theorem}

{\it Proof}. ${\rm (i)}\Longrightarrow {\rm (ii)}$: Let $M$ be an
arithmetical $R$-module and we show that the lattice of the
submodules of $M$ is distributive, i.e., for all submodules $K, L,
N$ of $M$, we have
$$(K+L)\cap N=(K+N)\cap(L+N).$$
To do this, in view of  \cite[Corollary 3.4 and Proposition
3.8]{AM}, it is enough for us to show that for all maximal ideals
$\frak m$ of $R$, we have
$$(K_{\frak m}+L_{\frak m})\cap N_{\frak m}=(K_{\frak m}+N_{\frak m})\cap(L_{\frak m}+N_{\frak m}).$$
Since the submodules of $M_{\frak m}$ are linearly ordered with
respect to inclusion, it follows that either $K_{\frak m}\subseteq
N_{\frak m}$ or $N_{\frak m}\subseteq K_{\frak m}$. Now, the desired
result easily follows from this and the modular law (see
\cite[Proposition 1.2]{SV}).

${\rm (ii)}\Longrightarrow {\rm (i)}$: Let $M$ be a distributive
$R$-module. Then, it easily follows from \cite[Ex. 9.11]{Sh} that,
for all maximal ideals $\frak m$ of $R$, the $R_{\frak m}$-module
$M_{\frak m}$ is also distributive. So without loss of generality we
may assume that $R$ is a quasi-local ring with the unique maximal
ideal $\frak m$, and we must show that the the submodules of $M$ are
linearly ordered with respect to inclusion. To do this, suppose that
the contrary is true, i.e., there exist two submodules $K$ and $L$
of $M$ such that $K\nsubseteq L$ and $L\nsubseteq K$. Then there
exist elements $x, y\in M$ such that $x\in K\backslash L$ and $y\in
L\backslash K$. Now, since $M$ is a distributive module, it follows
that
$$R(x+y)\cap Rx+R(x+y)\cap Ry=R(x+y)\cap (Rx+Ry)=R(x+y),$$
and so there are elements $w_1, w_2\in M$ such that
\begin{center}
$x+y=w_1+w_2$, where $w_1\in R(x+y)\cap Rx$ and $w_2\in R(x+y)\cap
Ry.$
\end{center}
Therefore, there exist elements $r, s, t\in R$ such that
$w_1=rx=s(x+y)$ and $w_2=ty$. Hence $(r-s)x=sy$, and so the elements
$r-s$ and $s$ are not units. This means that $r, s \in \frak m$. On
the other hand, since $x+y=rx+ty$, it follows that $(1-r)x=(t-1)y$,
and so $1-r\in \frak m$, which is a contradiction. Consequently one
of $K\subseteq L$ and $L\subseteq K$ must hold.

${\rm (i)}\Longrightarrow {\rm (iii)}$: Let $M$ be an arithmetical
$R$-module and let $K, L, N$ be submodules of $M$ such that $N$ is
finitely generated. We show that
$$(K+L):_RN=(K:_R N)+(L:_R N).$$
To do this, in view of  \cite[Corollary 3.4 and Proposition
3.8]{AM}, we may assume that $R$ is a quasi-local ring. Then the
submodules of $M$ are linearly ordered with respect to inclusion.
Hence, without loss of generality we may assume that $K\subseteq L$.
Then
\begin{center}
 $(K+L):_RN=L:_R N$ and $K:_R N\subseteq L:_R N.$
\end{center}
Now, the assertion follows.

${\rm (iii)}\Longrightarrow {\rm (i)}$: In view of  \cite[Corollary
3.4 and Proposition 3.8]{AM}, we may assume that $R$ is a
quasi-local ring with the unique maximal ideal $\frak m$. To
establish (i), suppose, on the contrary, that $M$ is not
arithmetical, and seek a contradiction. Then, there exist two
submodules $K$ and $L$ of $M$ such that $K\nsubseteq L$ and
$L\nsubseteq K$. Thus, there exist elements $x, y\in M$ such that
$x\in K\backslash L$ and $y\in L\backslash K$. By the hypothesis we
have
$$(Rx+Ry):_RR(x+y)=(Rx:_RR(x+y))+(Ry:_RR(x+y)).$$
hence $R=(Rx:_RR(x+y))+(Ry:_RR(x+y))$, and so $1=a+b$, where $a\in
(Rx:_RR(x+y))$ and $b\in (Ry:_RR(x+y))$. Therefore, there exist
elements $r, s\in R$ such that $ax+ay=rx$ and $bx+by=sy$. Hence
$(a-r)x=ay$ and $bx=(s-b)y$. Since $x\not\in Ry$ and $y\not\in Rx$,
it follows that $a, b \in \frak m$, so that $1\in \frak m$, which is
a contradiction.

${\rm (i)}\Longrightarrow {\rm (iv)}$: Let $M$ be an arithmetical
$R$-module and suppose that $K, L, N$ are  submodules of $M$ such
that $L$ and $N$ are finitely generated. It is clear that $K:_R
L\subseteq K:_R (L\cap N)$ and $K:_RN\subseteq K:_R(L\cap N),$ and
so $(K:_RL)+(K:_RN)\subseteq K:_R (L\cap N)$. Now, in order to show
the opposite inclusion, in view of \cite[Proposition 3.8]{AM}, it is
enough for us to show that, for all maximal ideals $\frak m$ of $R$,
$$(K:_R(L\cap N)/(K:_RL)+(K:_RN))_{\frak m}=0.$$
To do this, we have $$(K:_R(L\cap N))_{\frak m}\subseteq K_{\frak
m}:_{R_{\frak m}}(L\cap N)_{\frak m}=K_{\frak m}:_{R_{\frak
m}}(L_{\frak m}\cap N_{\frak m}),$$ and in view of  \cite[Corollary
3.4 and Proposition 3.8]{AM}, $$((K:_RL)+(K:_RN))_{\frak
m}=(K_{\frak m}:_{R_{\frak m}}L_{\frak m})+(K_{\frak m}:_{R_{\frak
m}}N_{\frak m}).$$ Now, since the submodules of $M_{\frak m}$ are
linearly ordered (with respect to inclusion), we may assume that
$L_{\frak m}\subseteq N_{\frak m}$. Then
\begin{center}
$(K:_R(L\cap N))_{\frak m}\subseteq (K_{\frak m}:_{R_{\frak
m}}L_{\frak m})$ and $((K:_RL)+(K:_RN))_{\frak m}=(K_{\frak
m}:_{R_{\frak m}}L_{\frak m})$
\end{center}
Therefore $$(K:_R(L\cap N))_{\frak m}\subseteq
((K:_RL)+(K:_RN))_{\frak m}.$$ Now the assertion follows easily from
\cite[Corollary 3.4]{AM}.

${\rm (iv)}\Longrightarrow {\rm (i)}$: According to the definition
we need to show that for every  maximal ideal $\frak m$ of $R$, the
submodules of the $R_{\frak m}$-module $M_{\frak m}$ are linearly
ordered (with respect to inclusion). To this end, suppose that the
contrary is true and look for a contradiction. Then, there exist two
submodules $K_{\frak m}$ and $L_{\frak m}$ of $M_{\frak m}$ (see
\cite[Ex. 9.11]{Sh}) such that  $K_{\frak m} \nsubseteq L_{\frak m}$
and $L_{\frak m} \nsubseteq K_{\frak m}$, where $K$ and $L$ are
submodules of $M$. Thus there are elements $x, y \in M$ such that
$x/1\in K_{\frak m}\backslash L_{\frak m}$ and $y/1\in L_{\frak
m}\backslash K_{\frak m}$. Now, by hypothesis (iv), we have
$$R=(Rx\cap Ry:_R Rx\cap Ry)=(Rx\cap Ry:_RRx)+(Rx\cap
Ry:_RRy)=(Ry:_RRx)+ (Rx:_RRy).$$ Hence $R_{\frak m}=(R_{\frak
m}y:_{R_{\frak m}}R_{\frak m}x)+ (R_{\frak m}x:_{R_{\frak
m}}R_{\frak m}y)$. But $R_{\frak m}y:_{R_{\frak m}}R_{\frak m}x$ and
$R_{\frak m}x:_{R_{\frak m}}R_{\frak m}y$ are proper ideals in
$R_{\frak m}$, we achieve a contradiction.

Finally, the equivalence between (ii) and (v) follows from
\cite[Proposition
7]{Ba}. \qed\\

The next result, which is the second main theorem of this section,
gives us two characterizations of strongly irreducible submodules of
the arithmetical modules in terms of primal and irreducible
submodules. To this end, recall that a proper submodule $N$ of an
$R$-module $M$ is called {\it primal submodule} of $M$ if
$Z_R(M/N)$, the set of zero-divisors of the $R$-module $M/N$, is an
ideal of $R$ (see \cite[Section 2, P. 193]{KN}). Then, it is easy to
see that $\frak p:=Z_R(M/N)$ is a prime ideal of $R$, called the
{\it adjoint prime ideal} of $N$. Also, in this case we say that $N$
is a $\frak p$-primal submodule of $M$.

\begin{theorem}
Let $M$ be an arithmetical $R$-module and let $N$ be a submodule of
$M$.  Then the following statements are equivalent:

${\rm (i)}$ $N$ is irreducible.

${\rm (ii)}$ $N$ is strongly irreducible.

${\rm (iii)}$ $N$ is primal.
\end{theorem}

{\it Proof}. ${\rm (i)}\Longrightarrow {\rm (ii)}$: Let $N$ be an
irreducible submodule of $M$ and suppose that $K$ and $L$ are two
submodules of $M$ such that $K\cap L\subseteq N$. Then, as $M$ is
arithmetical, it follows from Theorem 3.1 that
$$N=N+(K\cap L)=(N+K)\cap (N+L).$$
Since $N$ is irreducible, it yields that either $N=N+K$ or $N=N+L$,
and so either $K \subseteq N$ or $L \subseteq N$, as required.

The implication ${\rm (ii)}\Longrightarrow {\rm (i)}$ follows from
Lemma 2.2. In order to show ${\rm (ii)}\Longrightarrow {\rm (iii)}$,
suppose that $N$ is a strongly irreducible and let $a, b\in
Z_R(M/N)$. Then then there exist elements $x, y\in M\backslash N$
such that $ax, by\in N$, and so $Rx\cap Ry\nsubseteq N$. Hence,
there exists $z\in Rx\cap Ry$ such that $z\not\in N$. Thus,
$z=rx=sy$ for some elements $r, s\in R$, and so $$(a-b)z=rax-sby\in
N.$$ Therefore $(a-b)(z+N)=N$, i.e., $a-b\in Z_R(M/N)$. Moreover,
for every $c\in R$, we have $ac(x+N)=N$, and so $ac\in Z_R(M/N)$.
This shows that $Z_R(M/N)$ is an ideal of $R$, and hence $N$ is a
primal submodule.

 ${\rm (iii)}\Longrightarrow {\rm (ii)}$: Let $N$ be a primal submodule of
$M$. Then $Z_R(M/N)$ is a prime ideal of $R$; say $\frak
p:=Z_R(M/N)$. Suppose that $S$ is the multiplication closed subset
$R\setminus \frak p$ of $R$. It is then easily seen that $S(N)=N$.
Now, since the submodules of the $R_{\frak p}$-module $M_{\frak p}$
are linearly ordered with respect to inclusion, it follows that
$N_{\frak p}$ is a strongly irreducible submodule of $M_{\frak p}$.
Hence, in view of Lemma 2.5,
$N$ is also a strongly irreducible submodule of $M$, as required. \qed\\

\section{Strongly irreducible submodules in Noetherian modules}
The purpose of this section is to give a characterization for a
finitely generated module $M$ over a Noetherian ring $R$ to have a
strongly irreducible submodule. The main goal is Theorem 4.7. To
this end, as an application of Theorem 2.11, we first prove the
following proposition which is needed in the proof of that theorem.

\begin{proposition}
Let $R$ be a Noetherian ring and let $M$ be a finitely generated
$R$-module. Suppose that $N$ is a strongly irreducible submodule and
assume that $\frak pM_{\frak p}\neq N_{\frak p}$, where $\frak
p={\rm Rad}(N:_RM)$. Then the following conditions are hold:

${\rm (i)}$ The submodule $(N:_M\frak p)_{\frak p}$ of $M_{\frak p}$
is cyclic.

${\rm (ii)}$ $N_{\frak p}={\frak p}(N_{\frak p}:_{M_{\frak p}}\frak
pR_{\frak p})$.

${\rm (iii)}$ For each submodule $K$ of $M$ either $K\subseteq N$ or
$(N_{\frak p}:_{M_{\frak p}}\frak pR_{\frak p})\subseteq K_{\frak
p}$.
\end{proposition}

{\it Proof}. First of all, we note that as $N$ is a strongly
irreducible submodule of $M$, it follows from Lemma 2.2 that $N$ is
a $\frak p$-primary submodule of $M$. Moreover, in view of
Proposition 2.6, $N_{\frak p}$ is a strongly irreducible submodule
of $M_{\frak p}$. Also, since $\frak p={\rm Rad}(N:_RM)$ and $R$ is
Noetherian, it follows that there exists an integer $n\geq 1$ such
that $\frak p^nM \subseteq N$. Then, it is easy to see that
$N_{\frak p}\subsetneqq N_{\frak p}:_{M_{\frak p}}\frak pR_{\frak
p}$. Hence, it follows from Theorem 2.11 that $N_{\frak p}=\frak
p(N_{\frak p}:_{M_{\frak p}}\frak pR_{\frak p})$ and the submodule
$(N_{\frak p}:_{M_{\frak p}}\frak pR_{\frak p})$ of $M_{\frak p}$ is
cyclic.

Now, in order to show (iii), let $K$ be an arbitrary submodule of
$M$ such that $(N_{\frak p}:_{M_{\frak p}}\frak pR_{\frak
p})\nsubseteqq K_{\frak p}$. Then, it follows from Theorem 2.11 that
$K_{\frak p}\subseteq N_{\frak p}$, and hence as $N$ is $\frak
p$-primary it is easy to see that $K\subseteq N$. \qed \\

The next result of this section investigates whenever $N$ is a
strongly irreducible submodule of a module $M$ over a local
(Noetherian) ring $(R, \frak m)$, then $N$ and $N:_M\frak m$ are
comparable (under containment) to all submodules of $M$.
\begin{theorem}
Let $(R, \frak m)$ be a local (Noetherian) ring and let $M$ be a
finitely generated $R$-module. Suppose that $N$ is a strongly
irreducible submodule of $M$ such that $\frak mM\neq N$ and ${\rm
Rad}(N:_RM)=\frak m$. Then $N$ and $N:_M\frak m$ are comparable by
set inclusion to all submodules of $M$. In fact,
\begin{center}
$N=\bigcup \{K|\, K\, \text{is a submodule of}\,\, M\,
\text{and}\,\, K\subsetneqq N:_M\frak m\},$
\end{center}
 and
\begin{center}
$N:_M\frak m=\bigcap \{L|\, L\, \text{is a submodule of}\,\, M\,
\text{and}\,\, N\subsetneqq L\}.$
\end{center}
\end{theorem}
{\it Proof}. Let $K$ be an arbitrary submodule of $M$. We must show
that either $N\subseteq K$ or $K \subseteq N$. To do this, in view
of Theorem 2.11 either $K \subseteq N$ or $N:_M\frak m \subseteq K$.
 Now, if $K \nsubseteq N$,
then it follows that $N:_M\frak m \subseteq K$, and so $N\subseteq
K$. That is $N$ is comparable to all submodules of $M$. Also, if
$N:_M\frak m \nsubseteq K$, it follows that $K \subseteq N$, and
therefore $K\subseteq N:_M\frak m$. Thus $N:_M\frak m$ is also
comparable to all submodules of $M$. Now, we show that $N=\bigcup
\{K|\, K\, \text{is a submodule of}\,\, M\, \text{and}\,\,
K\subsetneqq N:_M\frak m\}.$ To do this, if $K$ is an arbitrary
submodule of $M$ such that $K\subsetneqq N:_M\frak m$, then
$N:_M\frak m\nsubseteqq K$, and so it follows from Theorem 2.11 that
$K \subseteq N$. Therefore
\begin{center}
$\bigcup \{K|\, K\, \text{is a submodule of}\,\, M\, \text{and}\,\,
K\subsetneqq N:_M\frak m\}\subseteq N.$
\end{center}

On the other hand, in view of Lemma 2.2 and the hypothesis ${\rm
Rad}(N:_RM)=\frak m$, we deduce that $N$ is an $\frak m$-primary
submodule of $M$, so that there exists an integer $k\geq 1$ such
that $\frak m^kM\subseteq N$. Hence we obtain that $N\subsetneqq
N:_M\frak m$, and so

\begin{center}
$N\subseteq \bigcup \{K|\, K\, \text{is a submodule of}\,\, M\,
\text{and}\,\, K\subsetneqq N:_M\frak m\}.$
\end{center}

Finally, in order to show
\begin{center}
$N:_M\frak m=\bigcap \{L|\, L\, \text{is a submodule of}\,\, M\,
\text{and}\,\, N\subsetneqq L\},$
\end{center}
if $L$ is an arbitrary submodule of $M$ such that $N\subsetneqq L$,
then in view of Theorem 2.11 we have $N:_M\frak m\subseteq L$. Hence
\begin{center}
$N:_M\frak m\subseteq\bigcap \{L|\, L\, \text{is a submodule of}\,\,
M\, \text{and}\,\, N\subsetneqq L\},$
\end{center}
Moreover, since $N\subsetneqq N:_M\frak m$ it follows that
\begin{center}
$\bigcap \{L|\, L\, \text{is a submodule of}\,\, M\, \text{and}\,\,
N\subsetneqq L\}\subseteq N:_M\frak m$,
\end{center}
and this completes the proof.  \qed\\

As an application, we derive the following consequence of Theorem
4.2, which shows that every strongly irreducible submodule $N$ of
finitely generated module $M$ over a Noetherian ring $R$ of
dimension one, is a distributive submodule, whenever the ideal
$N:_RM$ of $R$ contains a regular element on $M$. Recall that a
submodule $N$ of an $R$-module $M$ is called a {\it distributive
submodule} if for all submodules $L$ and $K$ of $M$, $(K\cap L)+
N=(K+ N)\cap(L+ N).$

\begin{corollary}
Let $R$ be a Noetherian ring and let $M$ be a finitely generated
$R$-module such that $\dim M=1$. Suppose that $N$ is a strongly
irreducible submodule of $M$ such that the ideal $N:_RM$ contains a
regular element on $M$. Then $N$ is a distributive submodule.
\end{corollary}
{\it Proof}. According to the definition it is enough to show that
for all submodules $L$ and $K$ of $M$, we have $$(K\cap L)+ N=(K+
N)\cap(L+ N).$$ To do this, it suffices to check the equation
locally at each prime ideal $\frak p$ in ${\rm Supp}(M)$. Now, if
$N:_RM\nsubseteq \frak p$, then $\frak p\not\in {\rm Supp}(M/N)$,
and so $N_{\frak p}=M_{\frak p}$. Hence the equation clearly holds
in this case. We therefore may assume that $N:_RM\subseteq \frak p$.
Then, it follows from $N:_RM\nsubseteq Z_R(M)$ that ${\rm
ht}_M(N:_RM)={\rm ht}_M\frak p=1$. Hence ${\rm Rad}(N:_RM)=\frak p$,
and therefore ${\rm Rad}(N_{\frak p}:_{R_{\frak p}}M_{\frak
p})=\frak pR_{\frak p}$. On the other hand, in view of Proposition
2.6 the submodule $N_{\frak p}$ is strongly irreducible in $M_{\frak
p}$. Consequently, without loss of generality we may assume that
$(R, \frak p)$ is local. Now, in view of Theorem 4.2, either
$K\subseteq N$ or $N\subseteq K$. If $K\subseteq N$, then
\begin{center}
$(K\cap L)+ N=N$ and $(K+ N)\cap(L+ N)=N\cap (L+N)=N,$
\end{center}
and so the equation holds in this case. Also, if $N\subseteq K$ then
by the modular law (see \cite[Proposition 1.2]{SV}), we have
$$(K\cap L)+ N=K\cap (L+N)=(K+N)\cap (L+N),$$ as required.  \qed\\

The following proposition gives us a characterization of strongly
irreducible submodule in a multiplication module over a commutative
Noetherian ring.

\begin{proposition}
Let $R$ be a Noetherian ring and let $M$ be a finitely generated
multiplication $R$-module. Suppose that $N$ is a proper submodule of
$M$ such that the ideal $N:_RM$ of $R$ is not prime. Then $N$ is
strongly irreducible if and only if there exists a submodule $L$ of
$M$ and a prime ideal $\frak p$ of $R$ such that $N\subsetneqq
L\subseteq \frak pM$ and that $N$ is $\frak p$-primary and for all
submodules $K$ of $M$ either $K\subseteq N$ or $L_{\frak p}\subseteq
K_{\frak p}$.
\end{proposition}
{\it Proof}. First, let $N$ be a strongly irreducible submodule of
$M$. Then it follows from Lemma 2.2 that $N$ is primary, and so
${\rm Rad}(N:_RM)$ is  a prime ideal of $R$. Let $\frak p={\rm
Rad}(N:_RM)$, and put $L=N:_M\frak p$. Now, in view of Proposition
4.1, it is enough for us to show that $N\subsetneqq L\subseteq \frak
pM$. To do this, as $R$ is Noetherian it follows that there exists
an integer $n\geq 1$ such that $\frak p^nM\subseteq N$. Let $t$ be
the least integer $n\geq1$ such that $\frak p^nM\subseteq N$. Thus
$\frak p^tM\subseteq N$ and $\frak p^{t-1}M\nsubseteq N$ (even if
$t=1$, simply because $N\neq M$), and so it yields that $N_{\frak
p}\subsetneqq N_{\frak p}:_{M_{\frak p}}\frak pR_{\frak p}$. Hence
$N\subsetneqq L$. In order to show $L\subseteq \frak pM$, since $M$
is a multiplication $R$-module, there exists an ideal $\frak a$ of
$R$ such that $L=N:_M\frak p=\frak a M$. Hence $\frak a \frak
pM\subseteq N$. Now, as the ideal $N:_RM$ is not prime, it follows
that $\frak pM\nsubseteq N$, and so $\frak a\subseteq \frak p$, note
that $N$ is $\frak p$-primary. Therefore $\frak aM\subseteq \frak
pM$, i.e., $L\subseteq \frak pM$.

In order to show the converse, let $L$ be a submodule of $M$ and
suppose that $\frak p$ is a prime ideal of $R$ such that
$N\subsetneqq L\subseteq \frak pM$ and assume that $N$ is $\frak
p$-primary. We show that $N$ is a strongly irreducible submodule. To
do this, let $T$ and $S$ be two submodules of $M$ such that $T\cap
S\subseteq N$. Now, suppose contrary is true, i.e., $T\nsubseteq N$
and $S\nsubseteq N$. Then, it follows from hypothesis that $L_{\frak
p}\subseteq T_{\frak p}\cap S_{\frak p}$. Whence $L_{\frak
p}\subseteq N_{\frak p}$, and so as $N$ is $\frak p$-primary it
easily follows that $L\subseteq N$, which is a
contradiction.  \qed\\

The following two propositions will serve to shorten the proof of
the main theorem of this section.

\begin{proposition}
Let $(R, \frak m)$ be a local (Noetherian) ring and let $M$ be a
non-zero finitely generated $R$-module such that ${\rm
deth}_R(R/\Ann_R(M))>0$ and that $\Gamma_{\frak m}(M)\neq0$. Suppose
$N$ is a strongly irreducible submodule of $M$ such that ${\rm
Rad}(N:_RM)=\frak m$. Then $\Gamma_{\frak m}(M)\nsubseteq N$.
\end{proposition}
{\it Proof}. In view of Proposition 4.1 the submodule $N:_M\frak m$
 of $M$ is cyclic, and so there exists $w\in M$ such that $N:_M\frak m=Rw$.
Now, in order to prove the claim, suppose the contrary is true, that
is $\Gamma_{\frak m}(M)\subseteq N$. Then $\Gamma_{\frak
m}(M)\subseteq Rw$, and so there exists an ideal $\frak a$ of $R$
such that $\Gamma_{\frak m}(M)= \frak aw$. On the other hand, since
$\Gamma_{\frak m}(M)$ is finitely generated it follows that there
exists $t\in \mathbb{N}$ such that $\frak m^t\frak aw=0$. Now, since
$N\subseteq Rw$ and the $R$-module $M/N$ has finite length, it
follows that the $R$-module $M/Rw$ has also finite length. Thus
there exists $l\in \mathbb{N}$ such that $\frak m^lM\subseteq Rw$.
Therefore $\frak m^{t+l}\frak aM\subseteq \frak m^t\frak aw=0$, and
so $\frak m^{t+l}\frak a\subseteq \Ann_R(M)$. Whence $$(\frak
m/\Ann_R(M))^{t+l}(\frak a+\Ann_R(M)/\Ann_R(M))=0,$$ and so $$(\frak
a+\Ann_R(M)/\Ann_R(M))\subseteq \Gamma_{\frak m}(R/\Ann_R(M)).$$
Now, since by hypotheses ${\rm deth}_R(R/\Ann_R(M))>0$, it follows
from \cite[Lemma 2.1.1]{BS} that $\Gamma_{\frak m}(R/\Ann_R(M))=0,$
and so $\frak a \subseteq \Ann_R(M)$. Hence $\frak aw=0$, and so
$\Gamma_{\frak
m}(M)=0$, which is a contradiction. \qed\\

\begin{proposition}
Let $R$ be a Noetherian ring and let $M$ be a finitely generated
$R$-module. Suppose $N$ is a strongly irreducible submodule of $M$
such that $\frak pM\nsubseteq N$ and that $\frak p \not\in
\Ass_R(R/\Ann_R(M))$, where $\frak p={\rm Rad}(N:_RM)$. Then, the
ideal  $N_{\frak p}:_{R_{\frak p}}M_{\frak p}$ of $R_{\frak p}$
contains a regular element on $M_{\frak p}$.
\end{proposition}
{\it Proof}. In view of Lemma 2.2 the submodule $N$ of $M$ is $\frak
p$-primary and hence it follows from $\frak pM\nsubseteq N$ that
$\frak pM_{\frak p}\nsubseteq N_{\frak p}$. So according to
Proposition 2.6 it may be assumed that $(R, \frak p)$ is a local
ring, and we must show that the ideal $\frak b:=(N:_RM)$ contains a
regular element on $M$. To this end, in view of \cite[Lemma
2.1.1]{BS}, it is enough for us to show that $\Gamma_{\frak
b}(M)=0$. Since ${\rm Rad}(\frak b)=\frak p$, it is suffices  to
establish that $\Gamma_{\frak p}(M)=0$.

In view of Proposition 4.1 either $\Gamma_{\frak p}(M)\subseteq N$
or $(N:_M\frak p)\subseteq \Gamma_{\frak p}(M)$. If $\Gamma_{\frak
p}(M)\subseteq N$, then the assertion follows from the proof of
Proposition 4.5.

 Therefore it may be assumed that $(N:_M\frak p)\subseteq \Gamma_{\frak p}(M)$.
 Then, there exists an integer $n\geq 1$ such that
$N\subseteq (0:_M\frak p^n)$, and so $\frak p^n N=0$. On the other
hand, since $\frak p={\rm Rad}(N:_RM)$, there is an integer $s\geq
1$ such that $\frak p^sM\subseteq N$. Hence $\frak p^{s+t}M\subseteq
\frak p^tN=0$, and so $\frak p^{s+t}M=0$. Therefore $M$ has finite
length. Now, it is easy to see that the $R$-module $R/\Ann_R(M)$ is
Artinian, and so $\Ass_R(R/\Ann_R(M))=\{\frak p\},$
which is a contradiction. \qed\\

Now we are prepared to prove the main theorem of this section, which
gives a characterization for a finitely generated module $M$ over a
Noetherian ring $R$ to have a strongly irreducible submodule.

\begin{theorem}
Let $R$ be a Noetherian ring and let $M$ be a finitely generated
$R$-module. Let $N$ be a submodule of $M$. Then, $N$ is strongly
irreducible if and only if $N$ is primary, $M_{\frak p}$ is an
arithmetical $R_{\frak p}$-module and $N=(\frak pM)^{(n)}$ for some
integer $n>1$, where $\frak p={\rm Rad}(N:_RM)$ such that $\frak
pM\nsubseteq N$ and $\frak p \not\in \Ass_R(R/\Ann_R(M))$.
\end{theorem}

{\it Proof}. First, let $N$ be a strongly irreducible submodule of
$M$. Then, it follows from Lemma 2.2 that $N$ is a primary submodule
of $M$. Now, let $\frak p:={\rm Rad}(N:_RM)$ and we show that
$M_{\frak p}$ is an arithmetical $R_{\frak p}$-module and $N=(\frak
pM)^{(n)}$ for some integer $n>1$. To do this end, in view of
Proposition 2.6 the submodule $N_{\frak p}$ of $M_{\frak p}$  is
strongly irreducible. Moreover, it is easy to see that $N=(\frak
pM)^{(n)}$ if and only if $N_{\frak p}=\frak p^nM_{\frak p}$.
Therefore, without loss of generality we may assume that $(R, \frak
p)$ is local, and that $(\frak pM)^{(n)}=\frak p^nM$. Then,
Proposition 4.6 shows that the ideal $(N:_RM)$ contains a regular
element on $M$, and according to Proposition 4.1 the submodule
$N:_M\frak p$ of $M$ is cyclic. Hence, there exists an element $w\in
M$ such that $N:_M\frak p=Rw$, and thus by Proposition 4.1 we have
$N=\frak pw$. Now, Since  by the Krull intersection theorem,
$\bigcap_{k\geq1}\frak p^kM=0$, and since $w\neq 0$, it follows that
$w\in \frak p^tM$, where $t$ is the greatest integer $i$ such that
$w\in \frak p^iM$. Then $N=\frak pw\subseteq \frak p^{t+1}M$. Now,
we show that $\frak p^tM$ is a cyclic submodule. To achieve this,
suppose the contrary is true. Then $\frak p^tM\neq\frak p^{t+1}M\cup
Rw$, and so there exists an element $y\in \frak p^tM$ such that
$y\not\in\frak p^{t+1}M\cup Rw$. Hence, it follows from $N\subseteq
\frak p^{t+1}M$ that $y\not\in N$. Also, one easily sees from
$y\not\in\frak p^{t+1} M\cup Rw$ that $w\not\in Ry$. Next, we show
that $Rw\cap Ry\subseteq N$. To do this, let $x\in Rw\cap Ry$. Then
there exist elements $a, b\in R$ such that $x=aw=by$. As $w\not\in
Ry$, it follows that $a\in \frak p$, and so $aw\in \frak pw$. This
shows that $x\in N$. Since $N$ is strongly irreducible, it yields
that either $Rw\subseteq N$ or $Ry\subseteq N$,  which is a
contradiction. Therefore $\frak p^tM$ is a cyclic submodule of $M$.
Hence, as $\frak p$ contains a regular element on $M$ (note that
$N:_RM$ contains a regular element on $M$ and $\frak p={\rm
Rad}(N:_RM)$), it follows from \cite[Theorem 2.3]{NB} that $R/{\rm
Ann}_R(M)$ is a ${\rm PID}$ and that $M\cong R/{\rm Ann}_R(M)$.
Consequently, every proper submodule of $M$ is of the form $\frak
p^vM$ for some integer $v\geq1$; and this shows that the submodules
of $M$ are linearly ordered with respect to inclusion, i.e., $M$ is
an arithmetical $R$-module and $N=\frak p^nM$ for some integer
$n>1$, as required.

In order to prove the converse, since $M_{\frak p}$ is an
arithmetical $R$-module, it follows that every submodule of
$M_{\frak p}$ is strongly irreducible. In particular $N_{\frak p}$
is strongly irreducible submodule of $M_{\frak p}$. Now, as $N$ is a
primary submodule of $M$, it follows from Lemma 2.5 that $N$ is
strongly irreducible.  \qed\\

\begin{corollary}
Let $R$ be a Noetherian integral domain and let $M$ be a non-zero
finitely generated torsion-free $R$-module. Then there exists a
strongly irreducible submodule $N$ of $M$ with $N:_RM$ is not prime
ideal of $R$ if and only if there is a prime ideal $\frak p$ of $R$
with $\frak pM\nsubseteq N$ and $M_{\frak p}$ is an arithmetical
$R_{\frak p}$-module.
\end{corollary}

{\it Proof}. Let $N$ be a strongly irreducible submodule of $M$ such
that $N:_RM$ is not prime ideal of $R$. Then, $N:_RM\neq 0$ and in
view of Lemma 2.2, $N$ is primary. Hence ${\rm Rad}(N:_RM)$ is a
prime ideal of $R$. Say $\frak p={\rm Rad}(N:_RM)$. Then, since $M$
is torsion-free, it is easy to see that $\frak p \not\in
\Ass_R(R/\Ann_R(M))$, and so in view of Theorem 4.7 the $R_{\frak
p}$-module $M_{\frak p}$ is arithmetical. Moreover, as the ideal
$N:_RM$ is not prime, one can easily check that $\frak pM\nsubseteq
N$.

Conversely, let $\frak p$ be a prime ideal of $R$ such that
$M_{\frak p}$ is an arithmetical $R_{\frak p}$-module. Then
$N:=(\frak pM)^{(2)}$ is a ${\frak p}$-primary submodule of $M$.
Whence in view of Theorem 4.7, $N$ is a strongly irreducible
submodule of $M$. Also, the ideal $N:_RM$ is not a prime ideal in
$R$, because if $N:_RM$ is a prime ideal, then $$N:_RM={\rm
Rad}(N:_RM)={\rm
Rad}((\frak pM)^{(2)}:_RM)=\frak p,$$ and so $\frak pM\subseteq N$, which is a contradiction. \qed\\

Before we state the final result of this paper, recall that the {\it
radical of} a submodule $N$ of an $R$-module $M$, denoted by ${\rm
rad}_M(N)$, is defined as the intersection of all prime submodules
containing $N$.

\begin{proposition}
Let $R$ be a Noetherian ring and let $M$ be a finitely generated
$R$-module. Let $N$ be a strongly irreducible submodule of $M$ such
that ${\rm rad}_M(N)=N$. Then $N$ is a prime submodule of $M$.
\end{proposition}
{\it Proof}. In view of Lemma 2.2 the submodule $N$ of $M$ is
primary and hence ${\rm Rad}(N:_RM)$ is a prime ideal of $R$. Hence
by applying \cite[Theorem 5]{Lu2} we deduce that $(N:_RM)$ is a
prime prime ideal of $R$, say $\frak p:=(N:_RM)$, and so in view of
\cite[Theorem 6.6]{Mat} we have ${\rm Ass}_R(M/N)=\{\frak p\}$.
Hence it follows from \cite[Corollary 9.36]{Sh} that $Z_R(M/N)=\frak
p$. Now, it is easy to see that $N$ is a prime submodule of $M$.
\qed\\
\begin{center}
{\bf Acknowledgments}
\end{center}
The authors are deeply grateful to Professors L. J. Ratliff, Jr.,
and R. Nekooei, for reading of the original manuscript and valuable
suggestions. Also, we would like to thank the Institute for Research
in Fundamental Sciences (IPM), for the financial support.


\end{document}